\documentclass[11pt]{amsart}
\usepackage{amsmath, amscd, amssymb}

\newcommand{\cA}{{\mathcal A}}

\newcommand{\cL}{{\mathcal L}}
\newcommand{\cO}{{\mathcal O}}
\newcommand{\cX}{{\mathcal X}}

\newcommand{\fh}{{\mathfrak h}}
\newcommand{\fk}{{\mathfrak k}}
\newcommand{\fm}{{\mathfrak m}}

\newcommand{\CC}{{\mathbb C}}
\newcommand{\FF}{{\mathbb F}}

\newcommand{\NN}{{\mathbb N}}
\newcommand{\QQ}{{\mathbb Q}}
\newcommand{\RR}{{\mathbb R}}
\newcommand{\ZZ}{{\mathbb Z}}

\newcommand{\fM}{{\mathfrak M}}
\newcommand{\fX}{{\mathfrak X}}
\newcommand{\fY}{{\mathfrak Y}}

\newcommand{\alg}{{\mathrm{alg}}}

\newcommand{\Cp}{{\CC_p}}
\newcommand{\Qp}{{\QQ_p}}

\newtheorem{theorem}{Theorem}[section]
\newtheorem{lemma}[theorem]{Lemma}

\newtheorem{lem}[theorem]{Lemma}
\newtheorem{proposition}[theorem]{Proposition}

\newtheorem{cor}[theorem]{Corollary}
\newtheorem{corollary}[theorem]{Corollary}
\newtheorem{prop}[theorem]{Proposition}
\newtheorem{conj}[theorem]{Conjecture}

\theoremstyle{remark}
\newtheorem{remark}[theorem]{Remark}

\theoremstyle{definition}
\newtheorem{definition}[theorem]{Definition} 

\begin{document}

\title[Local Andr\'{e}-Oort for the universal abelian variety]{Local Andr\'{e}-Oort conjecture for the universal abelian variety}
\author{Thomas Scanlon}
\thanks{Partially supported by NSF grant DMS-0301771 and an Alfred P. Sloan Fellowship.}
\address{University of California, Berkeley \\
Department of Mathematics \\
Evans Hall \\
Berkeley, CA 94720-3480 \\
USA}
\email{scanlon@math.berkeley.edu}
\maketitle
\begin{abstract}
We prove a $p$-adic analogue of the Andr\'{e}-Oort conjecture for 
subvarieties of the universal abelian varieties containing a dense set
of special points.  
Let $g$ and $n$ be integers with $n \geq 3$ and $p$ a prime number 
not dividing $n$.  Let $R$ be a finite extension of $W[{\mathbb F}_p^{\mathrm alg}]$, the 
ring of Witt vectors of the algebraic closure of the field of $p$ elements.
The moduli space $\cA = \cA_{g,1,n}$ of $g$-dimensional 
principally polarized abelian varieties with full level $n$-structure as
well as the universal abelian variety $\pi:\cX \to \cA$ over $\cA$ may 
be defined over $R$.  We call a point $\xi \in \cX(R)$ \emph{$R$-special}
if $\cX_{\pi(\xi)}$ is a canonical lift and $\xi$ is a torsion point of its
fibre.  We show that an irreducible  subvariety of $\cX_R$ containing a 
dense set of $p$-special points must be a special subvariety in the sense 
of mixed Shimura varieties.  Our proof employs the model theory of 
difference fields. 
\end{abstract}

\section{Introduction}
The Andr\'{e}-Oort conjecture in its usual form asserts that if an irreducible subvariety 
of a Shimura variety 
contains a Zariski dense set of special points then it must itself be a special 
variety, that is to say, a component of the image of a Shimura variety under a Hecke 
correspondence (see~\cite{Oort} and the appendix of~\cite{Andre}).
In initially proposing his version of this conjecture, Andr\'{e} relaxed the usual 
definition of ``Shimura variety'' allowing, for instance, abelian varieties and more generally
universal abelian varieties over moduli spaces as ambient varieties 
(see~\cite{Pink} for a general development of mixed Shimura varieties).
The conjecture even 
with the restriction to pure Shimura varieties, and even restricted 
to the case that the ambient Shimura variety is the moduli space of principally 
polarized abelian varieties of dimension $g$ with full level $n$ (which we shall always 
assume to be at least three) structure, $\cA_{g,1,n}$
(which we write as $\cA$ when the subscripts are understood),
is still 
open, though important special cases have been proven by many people~\cite{Andre curve, Edix curve, 
Edix hilb, EdYa}
and Moonen has proven the natural $p$-adic analogue of the Andr\'{e}-Oort conjecture
for $\cA$~\cite{Moonen}.
We formulate and prove the $p$-adic version of the Andr\'{e}-Oort 
conjecture for the universal abelian variety over $\cA$.  

The special points on and the special subvarieties of $\cA$ may be interpreted in terms of 
moduli:  a special point is nothing other than the moduli point of a CM-abelian 
variety (with the requisite additional structure) and a special subvariety is a submoduli 
variety, that is, it represents a submoduli problem given by some constraints on the 
Hodge structure of the abelian variety.   Likewise, the special points on and 
special subvarieties of the universal abelian variety $\pi:\cX \to \cA$ over $\cA$ 
may be described in a similar manner.  A point $x \in \cX(\CC)$ is special 
if and only if $\cX_{\pi(x)}$ is a CM-abelian variety and $x$ is a torsion point of 
$\cX_{\pi(x)}(\CC)$.  The special subvarieties are, essentially, translates of universal 
abelian varieties over special subvarieties of $\cA$.  More precisely, if $Y \subseteq \cX_\CC$
is special, then $\pi Y \subseteq \cA_\CC$ is a special subvariety of the moduli space and
for some natural number $M$, the variety $[M](Y) \subseteq \cX_\CC$ is an abelian subscheme of
$\cX_{\pi Y}$.  As with special subvarieties of Shimura varieties in the traditional 
sense, there are other characterizations of \emph{special subvariety of $\cX$} in terms of its
complex uniformization and preservation under sufficiently many Hecke correspondences.

The universal abelian variety $\pi:\cX_{g,1,n} \to \cA_{g,1,n}$ is actually over
$\ZZ[\frac{1}{n}]$ so that as long as $p$ is a prime not dividing $n$, we may regard
$\pi:\cX \to \cA$ via base change as being over $R := W[\FF_p^{\mathrm{alg}}]$.  
For Moonen, the $p$-adic special points on $\cA$ were precisely the $R$-points 
corresponding to canonical lifts (of ordinary abelian varieties).
  For us, the $p$-adic special points on $\cX$,
or as we prefer to say, the \emph{$R$-special points} are the $R$-points  of $\cX$
which are torsion on canonical lift fibres.  We will have more to say about canonical 
lifts in the next section.  For the moment, we merely note that the key to our proof 
is an interpretation of \emph{canonical lift} in terms of a difference equation.

The proof presented here uses the model theory of difference fields as developed 
by Chatzidakis and Hrushovski in~\cite{ChHr}.  Certainly, the logical formalism is 
intrinsic to neither the problem nor the proofs, but it provides a coherent
framework for the study of difference equations, especially where issues of uniformity 
intervene.  

J.~F.~Voloch introduced me to this problem in the Spring of 1998 and the general outline
of the solution presented here was sketched in conversations with him.  
I thank Z. Chatzidakis for pointing out that while quantifier elimination fails 
in ACFA, in some sense, it holds in the limit, and, in particular, in the form needed for
the proof presented here.  I thank the organizers of the December 2003 meeting 
on Special Points in Shimura Varieties at the Lorentz Center in Leiden for giving me the
opportunity to speak about this material and to learn about the Andr\'{e}-Oort conjecture
from the experts. 

\section{Preliminaries}

We start by fixing some conventions and notation.  Until further notice is
given, $p$ is a fixed prime number, $k := \FF_p^\alg$ is the 
algebraic closure of the field of $p$ elements, $W[k]$ is the ring of 
Witt vectors over $k$ and 
$\sigma:W[k] \to W[k]$ is the Witt-Frobenius.  That is, $\sigma$ is the unique
ring endomorphism of $R$ satisfying the congruence $\sigma(x) \equiv x^p \pmod{p}$
universally.  We fix an embedding $W[k] \hookrightarrow \CC$ of $R$ into the complex
numbers and an extension of $\sigma$ to an automorphism of $\CC$ (which we will
continue to call $\sigma$) so that $(\CC,+,\times,\sigma)$ is an $\aleph_1$-saturated
model of $\mathrm{ACFA}_0$ (see~\cite{ChHr} for the essentials of the theory of
difference closed fields of characteristic zero).  We fix as well a finite ramified 
extension $R$ of $W[k]$ of ramification index $e$.  As the results we aim to prove for 
$R$ follow from the corresponding results for a finite extension of $R$, we may, and 
do, assume that $\sigma$ preserves $R$.  Note that on $R$ the automorphism 
$\sigma$ still satisfies the congruence $\sigma(x) \equiv x^p \mod \fm_R$, though 
this congruence need not hold modulo $p$.    We normalize the $p$-adic valuation $v$
on the quotient field of $R$ by $v(p) = 1 \in \ZZ$ so that the minimal 
positive valuation of an element of $R$ is $\frac{1}{e}$.   

For $g \geq 1$ a positive integer and $n \geq 3$ the moduli space of principally 
polarized abelian schemes of relative dimension $g$ with full level $n$ structure 
is a fine moduli space $\cA_{g,1,n}$ over $\mathrm{Spec}(\ZZ[\frac{1}{n}])$ 
(see Theorem 7.9 and the following remark in~\cite{GIT}).
As such, there is a universal abelian variety $\pi:\cX_{g,1,n} \to \cA_{g,1,n}$
(also over $\mathrm{Spec}(\ZZ[\frac{1}{n}])$). In what follows, we write $\cX$ ($\cA$, 
respectively) for $\cX_{g,1,n}$ ($\cA_{g,1,n}$, respectively) except on the rare
occasions that the subscripts play some r\^{o}le.  Supposing in addition that 
$p$ does not divide $n$, then via base change to $R$, we may regard $\cA_R$ as 
a moduli space over $R$ and $\pi:\cX_R \to \cA_R$ as a universal abelian scheme 
over $R$.  Likewise, via base change over our embedding $R \hookrightarrow \CC$, 
we may regard $\pi:\cX_\CC \to \cA_\CC$ over $\CC$.  This embedding allows us to 
treat $\cX(R)$ as a subset of $\cX(\CC)$.

We recall the definition of \emph{canonical lift}.

\begin{definition}
Let $k'$ be a perfect field of characteristic $p$.  
We write the Witt-Frobenius on $W[k']$ as $\rho:W[k'] \to W[k']$.
If $A$ is an abelian scheme over $W[k']$ we say that $A$ is a 
\emph{canonical lift} if the special fibre $A_{k'}$
is ordinary and there is an isogeny $\psi:A \to A^\rho$ 
whose restriction to the special fibre is the geometric 
Frobenius morphism $F:A_{k'} \to A_{k'}^{(p)}$.
\end{definition}

There are various equivalent formulations of \emph{canonical lift}.  For instance,
in the case that $k'$ is a subfield of ${\mathbb F}_p^\alg$, 
$A$ being a canonical lift is equivalent to the restriction map 
$\mathrm{End}(A) \to \mathrm{End}(A_{k'})$ being an isomorphism.  
(See the appendix of~\cite{Messing}).
The only use we have for these equivalences is the observation that if $A/R$ is 
a canonical lift, then $A$ is a CM-abelian variety.

\begin{definition}
\label{specialpoint}
We define a point $\xi \in \cX(\CC)$ to be \emph{special} if the abelian 
variety $\cX_{\pi(\xi)}$ has complex multiplication and $\xi$ is a torsion 
point of $\cX_{\pi(\xi)}(\CC)$.  We say that
a point $\xi \in \cX(R)$ is \emph{$R$-special} if it is a torsion 
point in $\cX_{\pi(\xi)}(R)$ and the abelian scheme  $\cX_{\pi(\xi)}$ 
is a canonical lift.
\end{definition}

 Note that via our inclusion $\cX(R) \subseteq \cX(\CC)$
every $R$-special point is special.  Moreover, every special point is 
actually algebraic. 

\begin{remark}
It may make more sense to define a point $\xi \in \cX(R)$ to be $R$-special if
$\xi$ is a torsion point of $\cX_{\pi(\xi)}$ and $\cX_{\pi(\xi)}$ is a 
\emph{quasi-canonical lift}, namely, isogenous to the canonical lift of its special fibre. 
One might even drop the hypothesis that the special fibre be ordinary and simply define
an $R$-special point to be an $R$-valued point of $\cX$ which happens to be special.  
  We adopt Definition~\ref{specialpoint} for the technical reason that the
hypotheses of Moonen's theorem include the condition that there be a dense set of canonical lifts, rather
than quasi-canonical lifts.  It is very plausible that the methods employed in~\cite{Moonen} generalize
to the case of varieties containing a dense set of $R$-rational moduli points quasicanonical lifts or even
that this generalization follows formally from the published theorem.  However, at this writing I do not
see how to achieve this.  
\end{remark} 

As with pure Shimura varieties, the universal 
abelian variety $\cX_\CC$ admits a complex analytic uniformization
as a double coset space for a real algebraic group.  If 
$\fh_g$ is the Siegel upper half space of symmetric $g \times g$ complex
matrices with positive definite imaginary part, then the real symplectic
group $\mathrm{Sp}_{2g}(\RR)$ acts transitively and faithfully on $\fh_g$ via the 
formula $\left( \begin{matrix} A & B \\ C & D \end{matrix} \right) \cdot \tau := 
(A \tau + B) (C \tau + D)^{-1}$.  It bears remarking that the group of symplectic
similitudes, $\mathrm{CSp}_{2g}(\RR)$, acts on $\fh_g$ via the same formula, but, of course, 
the action is no longer faithful.  

Via the standard embedding $\rho:\mathrm{CSp}_{2g} \hookrightarrow \mathrm{GL}_{2g}$
we have an action 
of $\mathrm{CSp}_{2g}$ on ${\mathbb G}_a^{2g}$ and with respect to this action we may 
form the semidirect product $G := {\mathbb G}_a^{2g} \rtimes \mathrm{CSp}_{2g}$.  
The group $G(\RR)$ acts on $\RR^{2g} \times \fh_g$ via 
$(x;\gamma) \cdot (z;\tau) := (x + \rho(\gamma)(z); \gamma \cdot \tau)$.   
Real analytically, we may  identify $\RR^{2g} \times \fh_g$ with $\CC^g \times \fh_g$ via
$(x_1, \ldots, x_{2g}; \tau_1, \ldots, \tau_g) \mapsto (x_1 + x_{g+1} \tau_1, \ldots, x_g + x_{2g} \tau_g; \vec{\tau})$.
Via this identification, the action of $G(\RR)$ on $\CC^g \times \fh_g$ is via complex analytic 
functions.   If $\Gamma_n := \{ g \in \mathrm{Sp}_{2g}(\ZZ) ~|~ 
g \equiv \mathrm{Id}_{2g} \pmod{n} \}$, then one may identify $\cX(\CC)$ with 
$(\ZZ^{2g} \rtimes \Gamma_n) \backslash (\CC^g \times \fh_g)$.

It is sometimes convenient to 
work with an intermediate analytic covering.  Recalling the identification of $\RR^{2g} \times \fh_g$
with $\CC^g \times \fh_g$, we may see 
$\fX_g := \ZZ^{2g} \backslash (\CC^g \times \fh_g)$ as $(\RR/\ZZ)^{2g} \times \fh_g$
(again, real analytically) and the covering of $\cX(\CC)$ factors through
$\nu:\fX_g \to \cX$.  We write $\widetilde{\pi}:\fX_g \to \fh_g$ for the second projection 
function.  Note that $\widetilde{\pi}$ is holomorphic and descends to the 
map $\pi:\cX \to \cA$.

\begin{definition}
\label{specialsubvariety}
We say that an irreducible subvariety $Y \subseteq \cX_\CC$ is
\emph{special} if its image in $\cA_\CC$, $\pi Y$, is a special subvariety (or a
``variety of Hodge type'') in the sense of~\cite{Moonen} and for some positive integer $M$, the variety 
$[M](Y)$ is an abelian subscheme of $\cX_{\pi Y}$.
\end{definition}

\begin{remark}
We adopt Definition~\ref{specialsubvariety} because amongst the various equivalent definitions of
\emph{special subvariety} it fits most closely to our method of proof.  One might define 
a special subvariety to be a component of an image of a mixed Shimura variety under a generalized
Hecke correspondence (see~\cite{Pink aoml}).  What we call a \emph{special subvariety} would be 
a \emph{weakly special subvariety} containing a special point in Pink's terminology. Towards the
end of our proof, we implicitly use another characterization of special subvarieties, namely, that they
are preserved by Hecke correspondences of the form $\emph{T}_{\ell^m}$ (see below for a discussion 
of the Hecke correspondences on $\cX$) for arbitrarily large 
primes $\ell$ (and positive integers $m \in \ZZ_+$, possibly depending on $\ell$).  
\end{remark}

For any positive integer $m \in \ZZ_+$ the  group element $\gamma_m$
be the analytic automorphism of $\CC^g \times \fh_g$ given by 
$(z;\tau) \mapsto (z;\frac{1}{m} \tau)$
which descends to an algebraic correspondence $\widetilde{T}_m$ on 
$\cX$, and via $\pi:\cX \to \cA$, to an algebraic correspondence $T_m$ 
on $\cA$.  These correspondences have modular interpretations.   
Recall that via the principal polarization, the Weil pairing 
$\cX_x[m] \times \check{\cX}_x[m] \to \mu_m$ may be regarded as a 
pairing on $\cX_x[m]$ itself.  The correspondence $T_m$ may be described as 
follows.

\begin{eqnarray*}
\langle x, y \rangle \in T_m(\CC) & \Leftrightarrow & 
\begin{cases} \exists \psi:\cX_x \to \cX_y \\
	\text{ an isogeny of ppav's with level structure having } \\
   \cX_x[\psi] \leq \cX_x[m] \text{ a maximal
	isotropic subgroup of } \cX_x[m] \end{cases} 
\end{eqnarray*}

The isotropy hypothesis on $\cX_x[\psi]$ 
is equivalent to saying that $\# \cX_x[m] = m^g$ and the restriction of 
the Weil pairing to $\cX_x[\psi]$ is trivial.

The correspondence $\widetilde{T}_m$ may be expressed via the following
identity.

\begin{eqnarray*}
\langle x, y \rangle \in \widetilde{T}_m(\CC) 
& \Leftrightarrow & \begin{cases} \exists \psi:\cX_{\pi(x)} \to \cX_{\pi(y)} \\ \text{ an isogeny of ppav's with level structure having } \\

	\cX_{\pi(x)}[\psi] \leq \cX_{\pi(x)}[m] \text{ maximal isotropic and }
	\psi(x) = y \end{cases} 
\end{eqnarray*}

Note that if $A/R$ is a canonical lift and $\psi:A \to A^{\sigma^m}$ is an isogeny
lifting the Frobenius $F^m:A_k \to A^{(p^m)}_k$, then 
as the special fibre of $A[\psi]$ is $A_k[F^m]$, $A[\psi]$ is a 
maximal isotropic subgroup of $A[p^m]$ so that if $x \in \cA(R) \subseteq \cA(\CC)$
is a moduli point of $A$ (with some extra structure), then 
$\langle x, \sigma^m(x) \rangle \in T_{p^m}(\CC)$.

At this point we may state our main theorem.

\begin{theorem}
\label{main}
If $Y \subseteq \cX_\CC$ is an irreducible subvariety containing a Zariski 
dense set of $R$-special points, then $Y$ is a special subvariety of $\cX_\CC$.
\end{theorem}

As a special case of Theorem~\ref{main}, we have the following corollary about sections
of families of abelian varieties.

\begin{corollary}
\label{torsionsection}
Let $S$ be a variety over $R$ and $\pi:A \to S$ an abelian scheme over $S$
having the property that $A_s \not \cong A_t$ for $s \neq t \in S(R)$. 
Suppose that $f:S \to A$ is a section of $\pi$ and that the set of points $\xi \in S(R)$
for which $A_\xi$ is a canonical lift and $f(\xi) \in A_\xi(R)$ is torsion is Zariski 
dense in $S$.  Then $f$ is a torsion section.
\end{corollary}

\section{Proof of main theorem}
As mentioned in the introduction, we prove Theorem~\ref{main} by converting 
the problem into a question about difference varieties.  For $m \in \ZZ_+$ a 
positive integer we define

\begin{eqnarray*}
\Phi_m & := & \{ x \in \cX(\CC) ~|~ \langle x, \sigma^m(x) \rangle \in
	 \widetilde{T}_{p^m} (\CC) \} \\
\Psi_m & := & \{ x \in \cA(\CC) ~|~ \langle x, \sigma^m(x) \rangle \in
	 T_{p^m} (\CC) \} 
\end{eqnarray*}

\medskip

We prove Theorem~\ref{main} by analyzing sets of the form $Z(\CC) \cap \Phi_m$
where $Z \subseteq \cX_\CC$ is a subvariety.  The difference varieties $\Phi_m$
do not quite contain the $R$-special points, but they do not miss by much.

We recall a result on bounds on the number of torsion points near zero in an abelian 
scheme over a $p$-adic ring.   

\begin{proposition}
\label{mattuck}
Given a prime $\ell$ and a positive interger $r$ there is a 
bound $B = B(\ell,r) = \ell^{\lfloor \frac{r}{\ell -1} \rfloor}$
so that if $S$ is a discrete valuation 
ring of mixed characteristic $(0,\ell)$ of absolute ramification 
degree at most $e$ and $A$ is an abelian scheme over $S$, then 
every torsion point in the kernel of the reduction map 
$A(S) \to A_{\fk}(\fk)$ (where $\fk = S/\fm$ is the residue field)
has order dividing $B$.  Indeed, if $w$ is the discrete valuation on 
$S$ and $I := \{ x \in S ~|~ w(x) > \frac{w(\ell)}{\ell - 1} \}$, then 
reduction map $A(S) \to A(S/I)$ is injective on torsion.
\end{proposition}

\begin{proof}
Expressing the kernel of reduction in terms of the formal group law of $A$, this
is a slight elaboration of the theorems of Mattuck in~\cite{mattuck}.  
For the details, see Proposition 7 of~\cite{clark}.  
\end{proof}

With Proposition~\ref{mattuck} in place, we show that the $R$-special points 
satisfy a difference equation.

\begin{lemma}
\label{contains}
Let $m \in \ZZ_+$ be a positive integer and $B := B(p, e)$ be the 
bound of Theorem~\ref{mattuck}.  Set $\widehat{\Phi}_m := \{ x \in \cX(\CC) 
~|~ [B](x) \in \Phi_m \}$.  Then $\widehat{\Phi}_m$ contains all the $R$-special
points.
\end{lemma}

\begin{proof}
Let $\xi \in \cX(R)$ be an $R$-special point.  We write $X$ for $\cX_{\pi(\xi)}$.
Let $\psi:X \to X^\sigma$
be the lifting of the Frobenius $F:X_k \to X_k^{(p)}$
witnessing that $X$ is a canonical lift.   Then the isogeny
$\vartheta := \psi^{\sigma^{m-1}} \circ \cdots \circ \psi:X \to X^{\sigma^{m}}$
witnesses that $\pi(x) \in \Psi_m$. 

Now, on the $R$-points, both $\vartheta: \cX_{\pi(x)}(R) \to \cX_{\pi(x)}^{\sigma^m}(R)$
and $\sigma^m:\cX_{\pi(x)}(R) \to \cX_{\pi(x)}^{\sigma^m}(R)$ reduce to the Frobenius
map $F^m:X \to X^{(p^m)}$.  Thus, the image of the
group homomorphism $(\vartheta - \sigma^m): \cX_{\pi(x)}(R) \to \cX_{\pi(x)}^{\sigma^m}(R)$
is contained in the kernel of reduction.  Of course, this map takes torsion to torsion.
Hence, by our choice of $B$, we have $[B](\vartheta - \sigma^m)(\xi) = 0$.  That is, 
$\vartheta([B]\xi) = \sigma^m([B]\xi)$.  That is, $[B] \xi \in \Phi_m$ as claimed.
\end{proof}

While there are other points in $\widehat{\Phi}_m$, every algebraic point is special
(though not necessarily $R$-special).  For this proof we need a lemma on finite subgroups
of abelian varieties.

\begin{lem}
\label{isotropicrank}
Let $A$ be a principally polarized
abelian variety over $\CC$, $m > 1$ an integer, $\Xi \leq A[m](\CC)$ a 
maximal isotropic subgroup of the $m$ torsion on $A$, and $B \leq A$ an abelian 
subvariety.  Then $B(\CC) \cap \Xi$ is nontrivial.
\end{lem}
\begin{proof}
Let $\phi:\CC^g \to A(\CC)$ be a complex uniformization of $A(\CC)$ with kernel
$\Lambda$.  As $\Xi$ is maximal isotropic, via a change of variables of $\CC^g$, 
we may choose a symplectic basis $t_1, \ldots, t_g, \tau_1, \ldots, \tau_g$ 
for $\Lambda$ for which
 $t_i$ is the standard $i^\text{th}$ unit vector (for each $i \leq g$) 
and  $\phi^{-1} \Xi = \sum_{i=1}^g \ZZ t_i + \sum_{j=1}^g \ZZ \frac{\tau_j}{m}$.  
Let $E:\CC^g \to (\CC^\times)^g$ be the exponential map 
$(x_1, \ldots, x_g) \mapsto (\exp(2 \pi i x_1), \ldots, \exp(2 \pi i x_g))$.  
Then $\phi$ factors through $E$ giving a covering 
$\widetilde{\phi}:(\CC^\times)^g  \to A(\CC)$ for which $\widetilde{\phi}^{-1} \Xi$ is the 
group generated by $(\mu_m)^g$ and $E(\sum_{i=1}^g \ZZ \tau_i)$. 
The connected component $\widetilde{B}$ of $\widetilde{\phi}^{-1}(B(\CC))$ is an algebraic
torus of dimension $g' := \dim B$.   As such, $\widetilde{B} \cap \mu_m^g \cong 
(\ZZ/m\ZZ)^{g'}$.  As $\widetilde{\phi}$ is injective on torsion, 
 $\#(B(\CC) \cap \Xi) = m^{g'} > 1$ as claimed.
\end{proof}

\begin{lemma}
\label{algspecial}
For any positive integer $m$ the set $\widehat{\Phi}_m \cap \cX(\QQ_p^\alg)$ 
consists entirely of special points.
\end{lemma}

\begin{proof}
If $\xi \in \cX(\QQ_p^\alg)$, then for some $s \in \ZZ_+$ we have $\sigma^s(\xi) = \xi$.
Suppose in addition that $\xi \in \widehat{\Phi}_m$.  Let 
$\psi:\cX_{\pi(\xi)} \to \cX_{\pi(\xi)}^{\sigma^m}$ be an isogeny with 
$\psi([B]\xi) = \sigma^m([B]\xi)$ witnessing that 
$\xi \in \widehat{\Phi}_m$.  Set $\theta := \psi^{\sigma^{(s-1)m}} \circ \cdots \circ
\psi^{\sigma^m} \circ \psi: \cX_{\pi(\xi)} \to \cX_{\pi(\xi)}^{\sigma^{ms}}$.  Then 
$\theta$ witnesses that $\langle \pi(\xi), \pi(\xi) \rangle \in T_{sm}(\CC)$.
Any such point is a CM-moduli point.  Moreover, $[B](\theta - 1)(\xi) = 0$.
As the kernel of $\theta$ is a maximal isotropic
subgroup of $\cX_{\pi(\xi)}[p^{sm}]$, by Lemma~\ref{isotropicrank} the restriction of
$\theta$ to any abelian subvariety has positive degree.  In particular, 
the map $[B](\theta - 1)$ is an isogeny.  Thus, $\xi$ is a torsion point on a CM fibre.
That is, $\xi$ is a special point.
\end{proof}

Now that we know that the difference varieties $\widehat{\Phi}_m$ capture the 
$R$-special points, we analyze the structure of these difference varieties.
We begin with an observation about subgroups of abelian varieties defined 
by equations of the form $\psi(x) = \sigma^m(x)$ where $\psi$ is an 
isogeny.

\begin{lemma}
\label{locmod}
Let $m \in \ZZ_+$ be a positive integer, $A$ an abelian variety over $\CC$ and 
$\psi:A \to A^{\sigma^m}$ an isogeny.  We suppose that for every abelian subvariety
$B \leq A$ the restriction of $\psi$ to $B$ has a nontrivial kernel.  Then if 
$Y \subseteq A$ is an irreducible subvariety for which the set 
$\{ y \in Y(\CC) ~|~ \sigma^m(y) = \psi(y) \}$ is Zariski dense, then $Y$ is a 
translate of an abelian subvariety of $A$.
\end{lemma}

\begin{proof}
This follows from the main dichotomy theorem of~\cite{ChHr}. (See Section 4.1 of~\cite{Hr MM}.
For an algebraic proof, see Theorem 2.4 of~\cite{PiRo}.)
\end{proof}

We may restate Lemma~\ref{locmod} for general subvarieties of $A$.

\begin{cor}
\label{clm}
With $\psi:A \to A^{\sigma^m}$ as in Lemma~\ref{locmod}, if 
$Z \subseteq A$ is any subvariety, then the Zariski closure of 
the set $\{z \in Z(\CC) ~|~ \sigma^m(z) = \psi(z) \}$ is a finite 
union of translates of abelian subvarieties of $A$.
\end{cor}

We return now to the set $\widehat{\Phi}_m$.  For each point $a \in \Psi_m$, the
fibre $(\widehat{\Phi}_m)_a$ is a finite union of groups commensurable with groups of the form
$\{x \in \cX_a(\CC) ~|~ \sigma^m(x) = \psi(x) \}$ where $\psi$ ranges amongst 
the (finitely many) isogenies $\psi:\cX_a \to \cX_{\sigma^m(a)}$ with $\cX_a[\psi]$ a maximal 
isotropic subgroup of $\cX_a[p^m]$.   By Lemma~\ref{isotropicrank} all of these isogenies satisfy the hypothesis 
of Corollary~\ref{clm}.
Hence, combining this observation with 
Corollary~\ref{clm} we have the following proposition.

\begin{proposition}
\label{flm}
Let $m \in \ZZ_+$ be a positive integer and $Y \subseteq \cX_\CC$ any subvariety.
Then for any $a \in \cA(\CC)$ the set $Y(\CC) \cap (\widehat{\Phi}_m)_a$ is a
finite union of cosets of subgroups of $\cX_a(\CC)$.
\end{proposition}

We can do better than to say that the intersection $Y(\CC) \cap \widehat{\Phi}_m$
is fibrewise a finite union of cosets.  It follows on general principles 
that the degree (with respect to some fixed quasi-projective embedding) of the 
Zariski closure of $Y(\CC) \cap (\widehat{\Phi}_m)_a$ is bounded independently of 
$a \in \cA(\CC)$. It then follows that the abelian varieties appearing as the stabilizers
of the various components of the Zariski closure of $Y(\CC) \cap (\widehat{\Phi}_m)_a$
fall into finitely many algebraic families.  We approach the issue of describing the 
Zariski closure of $Y(\CC) \cap \widehat{\Phi}_m$ in a different manner.

We recall the definition of the Ueno locus.

\begin{definition}
Let $A$ be an abelian variety over $\CC$ and $Y \subseteq A$ a closed subvariety.
 We define the Ueno locus of $Y$, $\mathrm{Ueno}(Y)$, by 

\begin{eqnarray*}
\mathrm{Ueno}(Y)(\CC) & := & \begin{cases}  \{ y \in Y(\CC) ~|~ \exists B \leq A \\
	  	 \text{ an abelian subvariety for which } y + B \subseteq Y \} 
	\end{cases}
\end{eqnarray*}
\end{definition}

 From its definition, $\mathrm{Ueno}(Y)(\CC)$ is just a subset of the 
 $\CC$-rational points of $Y$.  However, it is actually a closed subvariety of 
  $Y$.  More is true: if $Y$ varies in an 
algebraic family, then so does $\mathrm{Ueno}(Y)$.

\begin{prop}
 \label{uniueno}
 If $Y \subseteq \cX_\CC$ is a closed subvariety of the universal abelian 
 variety, then there is a closed subvariety $\mathrm{Ueno}(Y) \subseteq Y$
 for which $\mathrm{Ueno}(Y)_a = \mathrm{Ueno}(Y_a)$ for all $a \in \cA(\CC)$.
 \end{prop}

\begin{proof}
 That the set 
 $\mathrm{Ueno}(Y)(\CC) := \{ y \in Y(\CC) ~|~ y \in \mathrm{Ueno}(Y_{\pi(y)}) \}$ is 
  constructible and that each fibre of $\mathrm{Ueno}(Y)$ is closed 
was noted in Lemma 11 of~\cite{Hr}.   We claim that
  $\mathrm{Ueno}(Y)$ itself is actually closed.  Write 
  $\mathrm{Ueno}(Y) = \bigcup Z_i \smallsetminus V_i$ where each $Z_i$ is a closed
  irreducible subvariety of $Y$ and $V_i \subsetneq Z_i$ is a proper closed subvariety.  
Let $Z = Z_i$ be one of the components of the Zariski closure of $\mathrm{Ueno}(Y)$.
Let $B$ be the fibrewise stabilizer of $Z$.  That is, 
$B(\CC) := \{x \in \cX(\CC) ~|~ \pi(x) \in \pi(Z) \& x + Z_{\pi(x)} = Z_{\pi(x)} \}$.
For a general point $a \in \pi(Z)(\CC)$ the fibre $B_a$ is a positive dimensional
algebraic subgroup of $\cX_a$.  As $B$ is a closed subvariety of $\cX$ and dimension
is uppersemicontinuous, it follows that for any $a \in \pi(Z)(\CC)$ the fibre
$B_a$ is a positive dimensional algebraic subgroup of $\cX_a$.  Thus,
for any such $a$ the fibre $Z_a$ is contained in the Ueno locus of $Y_a$.  
Thus, $\mathrm{Ueno}(Y)$ is closed.
   \end{proof}

So, our finiteness result on intersections of subvarieties of $\cX_\CC$ with 
$\widehat{\Phi}_m$ can take another form.

\begin{proposition}
\label{finiteoutueno}
If $m \in \ZZ_+$ is a positive integer and $Y \subseteq \cX_\CC$ is a closed subvariety, 
then for all $a \in \Psi_m$ the set 
$[Y(\CC) \smallsetminus \mathrm{Ueno}(Y)(\CC)] \cap (\widehat{\Phi}_m)_a$ is finite.
\end{proposition}

At this point, we interpret Proposition~\ref{finiteoutueno} in terms of algebraicity 
in difference fields.  Recall that we say that an element $b$ of some $\cL$-structure
$\fM$ (for some first-order language $\cL$) 
is model theoretically algebraic over the tuple $a = \langle a_1, \ldots, a_n \rangle$
if there is some $\cL$-formula $\phi(x;y_1, \ldots, y_n)$ having free variables 
$x, y_1, \ldots, y_n$ for which $\phi(b;a)$ holds in $\fM$ but there are only finitely 
many $b'$ from $\fM$ with $\phi(b';a)$ holding.  Proposition~\ref{finiteoutueno} says 
that the elements (or, more precisely, their coordinates) of
$[Y(\CC) \smallsetminus \mathrm{Ueno}(Y)(\CC)] \cap (\widehat{\Phi}_m)_a$ are
model theoretically algebraic over $a$ in $\cL(+,\times,\sigma^m,\vec{c})$ where 
$\vec{c}$ is a sequence of generators for a field of definition of $Y$.   
Recall (see Proposition 1.7 of~\cite{ChHr}) that in an existentially closed difference
field, the model theoretic algebraic closure of a set is the usual algebraic 
closure of the inversive difference field generated by that set.   
If $a \in \Psi_m$, then $\sigma^m(a)$ and $\sigma^{-m}(a)$ are algebraic over
$a$ (as $\langle a, \sigma^m(a) \rangle \in T_{p^m}(\CC)$ so that 
$\langle \sigma^{-m}(a), a \rangle \in T_{p^m}(\CC)$ as well).  Thus, in this 
case, being model theoretically algebraic over $a$ is the same as being 
algebraic in the usual sense.  The upshot of these observations is the following 
lemma.

\begin{lemma}
\label{uenofin}
Let $Y \subseteq \cX_\CC$ be a closed subvariety and $m \in \ZZ_+$ 
a positive integer.  Then there is a closed subvariety $Z \subseteq Y$
each component of which is generically finite over its image in $\cA_\CC$
for which $Y(\CC) \cap \widehat{\Phi}_m \subseteq \mathrm{Ueno}(Y)(\CC) 
\cup Z(\CC)$.
\end{lemma}

\begin{proof}
It suffices to work in a local affine chart $U \subseteq \cX_\CC$ and show that such a
$Z \subseteq Y \cap U$ exists.  Choose coordinates so that 
$U \subseteq {\mathbb A}^s \times {\mathbb A}^t$ with $\pi \upharpoonright U$ being 
given by the projection onto the first $s$ coordinates. By the above observations, for each 
$a \in \Psi_m \cap \pi U(\CC)$ there are polynomials $f_1(x_1,y_1, \ldots, y_s), 
\ldots, f_t(x_t,y_1, \ldots, y_s) \in \ZZ[x,y]$ so that $f_j(x;a) \not \equiv 0$ (for 
all $j \leq t$) and for any $b = \langle a; b_1, \ldots, b_t \rangle \in \widehat{\Phi}_m
\cap Y(\CC)$ we have $f_j(b_j;a) = 0$ (for all $j \leq t$).  Applying 
$\aleph_1$-compactness of $(\CC,\sigma)$ and taking products, we see that one may 
choose the polynomials $f_j$ independently from $a$.  Let $Z := Y \cap V(f_1(x_1,y),
\ldots, f_t(x_t,y))$.
\end{proof}

Lemma~\ref{uenofin}
permits us to presume that if $Y \subseteq \cX_\CC$ is an irreducible subvariety
containing a Zariski dense set of $R$-special points, then either $\pi \upharpoonright Y$
is generically finite  or $Y$ has a nontrivial Ueno fibration.  With the next lemma
we note that in the latter situation one may pass to a quotient thereby reducing to the 
former case.

\begin{lemma}
\label{quotient}
Suppose $Y \subseteq (\cX_{g,1,n})_\CC$ is a counter-example to Theorem~\ref{main}.
That is, $Y$ is a non-special irreducible subvariety containing a Zariski dense 
set of $R$-special points.  Then there is a counter-example $Z \subseteq (\cX_{g',1,n})_\CC$
with $\dim Z \leq \dim Y$ and $\pi \upharpoonright Z$ generically finite.
\end{lemma}

\begin{proof}
By Lemma~\ref{uenofin}, if this lemma is not already true with $Z = Y$, then 
$Y = \mathrm{Ueno}(Y)$.   As $\pi:\cX \to \cA$ is proper, $\pi Y \subseteq \cA$ is 
a closed subvariety.  Let $B := \mathrm{Stab}_{\cX_{\pi Y}}(Y)^0$ be the connected
component of the stabilizer of $Y$.  Then $B \leq \cX_{\pi Y}$ is an abelian subscheme 
of $\cX_{\pi Y}$ of relative dimension $d > 0$.    One would like to take $Z$ to 
be the quotient of $Y$ by $B$.  Unfortunately, while for every $a \in \pi Y(\CC)$ the
abelian variety $\cX_a$ is principally polarized, the 
quotient $\cX_a/B_a$ need not be principally polarized.
However, it does inherit a polarization from $\cX_a$ so that using Zahrin's 
trick (see, for instance, Remark 16.12 of~\cite{Milne}) 
we may regard it as an abelian subvariety of a principally polarized abelian variety 
of dimension $g' := 8 (g - d)$.  Let $\psi:Y \to (\cX_{g',1,n})_\CC$ be 
the map given by the quotient modulo $B$ followed by the inclusion given by 
Zahrin's trick.  Let $Z$ be the image of $Y$ under $\psi$.

As $Y$ contains a Zariski dense set of $R$-special points, it is defined over $R$.
Thus, $B$ is also defined over $Y$.  It follows that $\psi$ is defined over $R$. 
Thus, the $R$-special points in $Y$ are taken to $R$-special points in $Z$ and 
are therefore dense in $Z$.  Visibly, $Z$ has a fibrewise finite stabilizer.  
By Lemma~\ref{uenofin}, $\pi \upharpoonright Z$ is generically finite.
\end{proof}

So, the last interesting case to consider is that of $Y \subseteq \cX_\CC$
irreducible, containing a Zariski dense set of $R$-special points, and 
$\pi \upharpoonright Y$ generically finite.  With our next two lemmata we 
finish our analysis of this case and thereby conclude the proof of Theorem~\ref{main}.

Recall that for an algebraic variety $W$ over $\CC$ one may define the $\sigma$-Zariski
topology on $W(\CC)$ by taking the closed sets to be those defined locally by the 
vanishing of difference polynomials.  Recall that this topology is Noetherian. (See
Theorem 3.8.5 of~\cite{Cohn}.)

\begin{lemma}
If there is a counter-example $Y \subseteq \cX_\CC$ to Theorem~\ref{main},
with $\pi \upharpoonright Y$ generically finite, then there is 
such a counter-example with $Y(\CC) \cap \widehat{\Phi}_1$ irreducible in the 
$\sigma$-Zariski topology.
\end{lemma}

\begin{proof}
Let $Y \subseteq \cX_\CC$ be a potential counter-example to Theorem~\ref{main}. 
As the $\sigma$-Zariski topology on $Y$ is Noetherian, we may write 
$Y(\CC) \cap \widehat{\Phi}_1$ as a finite union $\bigcup_{i=1}^m Y_i$ of $\sigma$-irreducible
difference varieties.   For $N$ sufficiently large, we may write each $Y_i$ as
$$Y_i := \{ x \in \cX_\CC(\CC) ~|~ \langle x, \sigma(x), \ldots, \sigma^N(x) \rangle \in Z_i(\CC) \}$$ 
where $Z_i \subseteq \cX_\CC \times \cX_\CC^\sigma \times \cdots \times 
\cX_\CC^{\sigma^{N}}$ is an irreducible variety.  Passing to a subvariety
if need be, we may assume that 
$\{ \langle x, \sigma(x), \ldots, \sigma^N(x) \rangle ~|~ x \in Y_i \}$ is 
Zariski dense in $Z_i$.   

Since $\cX$ is defined over the fixed field of $\sigma$, we see that $\cX^{\sigma^i} = \cX$ 
for any $i$.  Moreover, by identifying the $(N+1)$-tuple of abelian varieties
encoded by an element of $\cX^{N+1}$ with their product, we 
may regard $(\cX_{g,1,n})^{N+1}$ as a subscheme of $\cX_{(N+1)g,1,n}$.  

As for each $a \in \widehat{\Phi}_1$ the elements
$\sigma^i(a)$ are algebraic over $a$, it follows that via the natural projection 
map $\rho: \cX^{N+1} \to \cX$ that each $Z_i$ is generically finite over $Y$.  
Moreover, because $Y$ is irreducible, the map $\rho \upharpoonright Z_i:Z_i \to Y$
is dominant.  

The function $x \mapsto \langle x, \sigma(x), \ldots, \sigma^N(x) \rangle$ 
takes $R$-special points to $R$-special points.  It follows that for at least 
one $i \leq m$ that the $R$-special points are dense in $Z_i$.  Fix one such and
set $Z := Z_i$.

Now, if $Z$ were special, then $Y$ would be special as well for $\pi Y$ would be
the image of the special variety $\pi Z$ under a map of Shimura varieties and 
on fibres $\rho$ is a map of group varieties so that it would take the torsion
variety $Z$ to a torsion variety.  Thus, if $Y$ is a counter-example
to Theorem~\ref{main}, so is $Z$.
\end{proof}

The key observation in the next lemma was pointed out to the author by Z. Chatzidakis.

\begin{lemma}
\label{torsionfibre}
If $Y \subseteq \cX_\CC$ is a subvariety containing a Zariski dense 
set of $R$-special points, for which the restriction of $\pi$ to $Y$
is generically finite, and $Y(\CC) \cap \widehat{\Phi}_1$ is $\sigma$-irreducible,
then there is a Zariski open and dense set $U \subseteq \pi Y$ so that
for every $R$-special point $\xi \in U(\CC)$ the fibre $Y_\xi(\CC)$ 
consists entirely of torsion points.
\end{lemma}

\begin{proof}
Let $d := \deg(\pi \upharpoonright Y)$ be the generic degree of $\pi$
restricted to $Y$ and set $m := d!$.  Let $K \subseteq \CC$ be a countable
algebraically closed inversive difference field over which $Y$ is defined.
As $Y$ contains a dense set of algebraic points, we could take $K = \QQ^\alg$. 
Let $a \in Y(\CC) \cap \widehat{\Phi}_1$
be a Weil-generic point.  That is, the transcendence degree of $K(a)$ over $K$ is 
equal to $\dim Y$.  Such a point exists by $\aleph_1$-saturation.  Set 
$a' := \pi(a) \in (\pi Y)(\CC) \cap \Psi_1$.  Let $b \in (\pi Y)(\CC) \cap \Psi_1$
be any generic point and $c \in Y_b(\CC)$ any point in the fibre above 
$b$ in $Y$.  It suffices to show that $c \in \widehat{\Phi}_m$.  Indeed, the 
set $\{ x \in (\pi Y)(\CC) \cap \Psi_1 ~|~ Y_x(\CC) \subseteq \widehat{\Phi}_m \}$
is definable.  So, if it contains all generic points, then it contains a 
(relatively) Zariski open and dense set.

As $Y(\CC) \cap \widehat{\Phi}_1$ is $\sigma$-irreducible, so is $(\pi Y)(\CC) \cap \Psi_1$.
Thus, the difference fields $L := K \langle a' \rangle_\sigma$ and 
$L' := K \langle b \rangle_\sigma$
are isomorphic over $K$ via an isomorphism $\rho$ taking $a'$ to $b$.  
Now, as $Y$ is an irreducible variety, the \emph{fields} $L(a)$
and $L'(c)$ are isomorphic over $\rho$ via $a \mapsto c$. 
We extend this map to 
$L(a)^\mathrm{alg} = L^\mathrm{alg}$
and call the extension $\widetilde{\rho}$.  

As $[L(a) : L] = d$, the 
normalization $M$ of this extension has degree dividing $d! = m$ over 
$L$.  Of course, the normalization $M'$ of 
$L'(c)$ over $L'$ is isomorphic
to $M$ via $\widetilde{\rho}$.  The group generated by 
$\sigma$ in $\mathrm{Aut}(L^\mathrm{alg})$ (respectively, $\mathrm{Aut}((L')^\mathrm{alg})$)
acts on $\mathrm{Gal}(M/L)$ (respectively, $\mathrm{Gal}(M'/L')$) by conjugation.
As the order of $\mathrm{Gal}(M/L) (respectively, \mathrm{Gal}(M'/L')$)
divides $m$, we see that 
$\sigma^m$ acts trivially.  In particular, on $M$, $\widetilde{\rho}$ and 
$\sigma^m$ commute.  

Now, $a \in \widehat{\Phi}_m$ meaning that 
$\langle [B] a, [B] \sigma^m (a) \rangle \in \widetilde{T}_{p^m}(\CC)$.  So that,
 
\begin{eqnarray*}
\widetilde{T}_{p^m}(\CC) & \ni & \widetilde{\rho} (\langle [B] a,  [B] \sigma^m(a) \rangle) \\
	& = & \langle [B] \widetilde{\rho}(a), [B] \widetilde{\rho} (\sigma^m(a)) \rangle \\
	& = & \langle [B] \widetilde{\rho}(a), [B] \sigma^m (\widetilde{\rho}(a)) \rangle \\
	& = & \langle [B] c, [B] \sigma^m(c) \rangle 
\end{eqnarray*}
That is, $c \in \widehat{\Phi}_m$ as claimed.

Therefore, for every generic $b \in (\pi Y)(\CC) \cap \Psi_1$ we have 
$Y_b(\CC) \subseteq \widehat{\Phi}_m$.  By compactness, there is a Zariski
dense and open  $U \subseteq \pi Y$ so that $Y_b(\CC) \subseteq \widehat{\Phi}_m$
for all $b \in U(\CC) \cap \Psi_1$.  By Lemma~\ref{algspecial}, if 
$\xi \in U(\CC) \cap \Psi_1$ is special, then every algebraic point in
$(\widehat{\Phi}_m)_\xi$ is special.  As $Y_\xi(\CC) \subseteq (\widehat{\Phi}_m)_\xi$,
we have that $Y_\xi(\CC)$ consists entirely of torsion points.
\end{proof}

With all of the above ingredients in place, we conclude the proof Theorem~\ref{main}.

\begin{proof}
By Lemma~\ref{torsionfibre} we may reduce to the case that $Y \subseteq \cX$ is 
an irreducible subvariety containing a dense set of $R$-special points,
with the restriction of $\pi$ to $Y$ being generically finite and for which there
is a Zariski open and dense set $U \subseteq \pi Y$ having the property that 
$Y_a(\CC)$ consists entirely of torsion points for all special points 
$a \in (\pi Y)(\CC) \cap \Psi_1$.  

As $\pi:\cX \to \cA$ is proper, the constructible set $\pi Y$ is a closed 
subvariety of $\cA$.  As the $R$-special points are dense in $Y$, they are 
also dense in $\pi Y$.  By the main theorem of~\cite{Moonen}, the variety $\pi Y$
is a special subvariety of $\cA$.   As $\pi Y$ is special, for arbitrarily 
large primes $\ell$ there are discrete valuation rings $S \subseteq \QQ^\mathrm{alg}$
of residue characteristic $\ell$ with 
the set 
$$C_{S,\ell,p} := \{ \zeta \in U(S) \cap U(R) ~|~ \zeta \text{ 
is canonical at $p$ and at $\ell$} \}$$
being Zariski dense in $\pi(Y)$ (see
Corollary 3.8 of~\cite{Moonen}).   If $\ell$ is 
sufficiently large, then $\pi \upharpoonright Y: Y \to \pi Y$ is (generically) 
unramified at $\ell$, say, over the Zariski open and dense set $V \subseteq \pi Y$.
Hence, for $\zeta \in C_{S,\ell,p} \cap V(S)$ the fibre $Y_\zeta(\CC)$
consists entirely of $\ell$-unramified torsion points.  Therefore, $Y$ contains 
a Zariski dense set of $S$-special points and is, thus, fixed by the Hecke
correspondence $\widetilde{T}_{\ell^m}$ for some positive integer $m$.  If 
$\fY \subseteq \fX_g$ is a component of $\nu^{-1}(Y(\CC))$, then $\epsilon := \gamma_{\ell^m} \delta$
fixes $\fY$ setwise for some $\delta \in \Gamma_n$.  Taking $\ell$ sufficiently large, 
by Theorem 6.1 of~\cite{EdYa} each orbit of $\epsilon$ on $\widetilde{\pi}(\fY)$
is dense.  Consequently, $\fY$ is a homogenous space for its stabilizer.  

Let $y \in Y(\CC)$ be a special point and $\tilde{y} \in \fY$ a lifting.  
Let $z \in Y(\CC)$ be any other point and $\tilde{z} \in \fY$ a lifting.
Using the real analytic trivialization of $\fX_g$, write $\tilde{y} = (x_y;\tau_y)$ and
$\tilde{z} = (x_z;\tau_z)$.  As $y$ is special, $x_y$ is a torsion point, say, of 
order $N \in\ZZ_+$. As $\fY$ is 
a homogenous space, there is some $\gamma$ in its stabilizer with $\gamma \tilde{y} = \tilde{z}$.
Fibrewise, the action on $(\RR/\ZZ)^{2g}$ is via group homomorphisms.  Hence, $x_z$ is also an 
$N$-torsion.  Therefore, $Y$ is a component of $\cX_{\pi Y}[N]$ and is, in particular,
a special variety.
\end{proof}

\section{Refinements of the method}
Our method of proof for Theorem~\ref{main} generalizes in several ways.
Firstly, using the model theory of valued difference fields, one can transfer 
Theorem~\ref{main} to include the case of $R$-special points where the residue field of 
$R$ is an arbitrary algebraically closed field of characteristic $p$.  
Secondly, the general theory
of difference equations produces explicit (albeit enormous) bounds on the degree of the Zariski
closure of solutions to difference equations in terms of geometric data describing 
the equations.  
Thirdly,
these methods apply to more general mixed Shimura varieties,
universal semi-abelian varieties, for instance. 

We explain these generalizations in less detail as we proceed. 
That is, we give complete statements and proofs for the generalization to arbitrary residue fields.  In the case of the explicit degree bounds, we 
compute the bounds relative to certain geometric degree computations.
We confine ourselves to the following comment about the more general mixed Shimura varieties.  
While it is clear that these methods are germane to the study of 
$p$-adic special points on the universal semiabelian variety, we lack a proven 
description of the subvarieties of the base moduli space which contain a 
dense set of $p$-adic special points.  

\subsection{Transcendental special points}
For the time being, $k'$ is an algebraically closed field of characteristic $p$ 
and $R'$ is a finite extension (of ramification index $e$) of the Witt vecotrs $W[k']$.
We let $\sigma:W[k'] \to W[k']$ be the relative Frobenius on the Witt vectors and 
we assume that $\sigma$ extends to an automorphism (which we continue to denote by $\sigma$)
of $R$.  We may, and do, express $R'$ as $W[k'][x]/(P)$ where $P \in \ZZ[x]$ is a polynomial
over the integers.    As before, let $k = \FF_p^\alg$ and set $R := W[k][x]/(P)$.  We write 
$K$ for the field of fractions of $R$ and $K'$ for the field of 
fractions of $R'$.  We write $v'$ for the valuation on $K'$ and $v$ for the valuation on $K$.
As before, these valuations are normalized so that $v(p) = v'(p) = 1$.

Definition~\ref{specialpoint} generalizes immediately to a give a notion of $R'$-special point.
A point $\xi \in \cX(R')$ is \emph{$R'$-special} if $\cX_{\pi(\xi)}$ is a canonical lift
and $\xi$ is a torsion point of $\cX_{\pi(\xi)}(R)$.    Definition~\ref{specialsubvariety} 
also generalizes immediately to $R'$, but the analogue of Theorem~\ref{main} would be false with this
immediate generalization.  Rather, we say that an irreducible subvariety
 $Y \subseteq \cX_{K'}$ is 
$R'$-special if $Y(K')$ contains an $R'$-special point, 
$\pi Y \subseteq \cA_{K'}$ is a totally geodesic subvariety in the sense of~\cite{Moonen}
and $[M](Y)$ is an abelian subscheme of $\cX_{\pi Y}$ for some positive integer $M$.
Equivalently, 
an $R'$-special variety is a weakly special subvariety in the sense of~\cite{Pink aoml} which contains
an $R'$-special point.  

With these definitions in place we may state the $R'$-version of Theorem~\ref{main}.

\begin{theorem}
\label{transmain}
If $Y \subseteq \cX_{K'}$ is an irreducible subvariety of $\cX_{K'}$ containing a Zariski dense set of 
$R'$-special points, then $Y$ is an $R'$-special subvariety.
\end{theorem}

\begin{remark}
Proposition 6.4 of~\cite{Moonen} may be seen as the special case of Theorem~\ref{transmain} where 
$Y$ is defined over $K$ and is equal to $\pi^{-1} (\pi Y)$.
\end{remark}

The key step in our proof of Theorem~\ref{transmain} is the fact that the extension of valued difference fields,
$(K,+,\times,v,\sigma) \hookrightarrow (K',+,\times,v',\sigma)$ is elementary in the sense of 
mathematical logic (see Corollary 10.4 of~\cite{Sc} or~\cite{BMS}).  On the face of it, Theorem~\ref{main}
is not first-order expressible, but we shall show that a failure of Theorem~\ref{transmain} may be converted
to a failure of Theorem~\ref{main}.

\begin{lemma}
\label{definable}
There is a formula $\phi$ in the language of valued difference fields for which 
the set of realizations of $\phi$ in $W[k'][\frac{1}{p}]$ is exactly the set of $R'$-special points in $\cX(K')$. [NB: This includes the case of $k' = k$.]
\end{lemma}

\begin{proof}
Strictly speaking, first-order formulae do not take points on quasi-projective schemes as 
arguments, but through the usual encoding tricks we may regard $\cX(R')$ together with the 
reduction  map $\cX(R') \to \cX(k')$ as part
 of the first-order structure.

More precisely, the ring $W[k']$ is definable in $W[k'][\frac{1}{p}]$ as the set of elements
of nonnegative valuation.  Using the presentation $R'$ as $W[k'][\frac{1}{p}][x]/(P)$,
we may regard any element of $R'$ as a $d = \deg(P)$-tuple of elements of $W[k']$
and the product and sum operations
on $R'$ may be expressed in terms of these coordinates.  
Moreover, one may indentify those elements of 
$R'$ which actually lie in $W[k]$
 as those $d$-tuples whose last $d-1$ coordinates are all zero.
The set of $R'$-rational points of $N$-dimensional
projective space may be expressed as the set of $(N+1)$-tuples of elements of $R'$
not all of 
which have positive valuation modulo the definable 
equivalence relation $x \sim y \Leftrightarrow (\exists \lambda) 
[v(\lambda) = 0 \ \& \ \lambda x = y]$.  
Using the above interpretation of $R'$ in $W[k']$,
this may be understood as a definable set of $d(N+1)$-tuples of elements of $W[k']$
modulo a definable 
equivalence relation.
Fixing a presentation of $\cX$ as a quasiprojective scheme over $\ZZ[\frac{1}{n}]$,
we may express the condition that
(the equivalence class of the $d(N+1)$-tuple of elements of $W[k]$ representing)
$y$ is an element of $\cX(R')$
in terms of the equations and inequations defining $\cX$.  Likewise, the 
reduction map may be expressed in terms of the standard interpretable reduction map on $R'$.
With these niceties in place, we work with variables ranging over $\cX(R')$ and speak about
the reduction map as if they are part of the first-order structure.

Let $\cA'$ be the moduli space of $g$-dimensional principally polarized abelian varieties with 
full level $n$ structure and a choice of a maximal isotropic subgroup of the $p$-torsion.  
There is a map $\vartheta:\cA' \to \cA$ which takes the moduli point of
the quadruple $(A,\lambda,\alpha,\Gamma)$
consisting of an abelian variety $A$, a polarization $\lambda$, level $n$-structure $\alpha$, and the maximal 
isotropic group of $p$-torsion $\Gamma$ to the moduli point of $A/\Gamma$ with the induced polarization 
and level structure.  Over $\QQ$, these are fine moduli spaces and the map $\vartheta$ lifts to a map on the universal 
abelian varieties $\widetilde{\vartheta}:(\cX')_\QQ \to \cX_\QQ$ which on the fibre over $[(A,\lambda,\alpha,\Gamma)]$ is
the quotient homomorphism $A \to A/\Gamma$. 

Let $\varpi:\cX'_\QQ \to \cX_\QQ$ be the natural projection map which forgets the choice of a $p$-torsion group.
If $\varpi(x) = y$, then we can consider $x$ as $(y,\Gamma)$ where $\Gamma$ is the $p$-torsion group 
encoded by $x$.  We write this $\Gamma$ as $\Gamma_x$.

Let $B = B(p,e)$ be the bound of Lemma~\ref{mattuck}.
Let $\phi(x)$ be the formula which asserts that $\pi(x) \in \cA(W[k'])$ and that there exists 
a point $\widetilde{x} \in \cX'(K')$ so that 
\begin{itemize}
\item $\varpi(\widetilde{x}) = x$
\item $\sigma([B]x) = \widetilde{\vartheta}([B]\widetilde{x})$
\item $\Gamma_{\widetilde{x}} \leq \cX_{\pi(x)}(W[k'])$ and the reduction map on $\Gamma_{\widetilde{x}}$
	is injective.
\end{itemize}

If $\phi(a)$ holds as witnessed by $\widetilde{x} = (a,\Gamma)$, then by the third condition the map 
$\widetilde{\vartheta}_{\pi(a)}:\cX_{\pi(a)} \to \cX_{\pi(\vartheta(\pi(a))}$ lifts the Frobenius and
the special fibre is ordinary. 
By the second condition, the codomain of this map is actually $\cX_{\pi(a)}^\sigma$.  Thus, 
$\cX_{\pi(a)}$ is a canonical lift.  Exactly as in the proof of Lemma~\ref{contains}, $a$ is 
a torsion point.  Hence, $a$ is an $R'$-special point.   Conversely, if $a$ is an $R'$-special
point witnessed by $\psi:\cX_{\pi(a)} \to \cX_{\pi(a)}^{\sigma}$ a lifting of the Frobenius, then as
the reduction map is injective on $[B] (\cX_{\pi(a)}(R')_{\mathrm{tor}})$ and $\psi$ and $\sigma$
agree upon reduction, $\sigma(a) = \vartheta(a,\ker \psi)$.

\end{proof}

With the next lemma we transfer Theorem~\ref{main} from $R$ to $R'$
for varieties defined over $K$.

\begin{lemma}
\label{firsttransfer}
Let $Y \subseteq \cX_K$ be an irreducible subvariety of $\cX_K$.  If $Y(K')$ contains 
a Zariski dense set of $R'$-special points, then $Y$ is a special subvariety of $\cX_K$.
\end{lemma}
\begin{proof}
If $Y$ is not a special subvariety, then by Theorem~\ref{main}, there is a 
proper closed subvariety $Z \subset Y$ (defined over $K$) so that every 
$R$-special point in $Y$ actually lies on $Z$.  By Lemma~\ref{definable}, 
this is expressible as a first-order sentence in the language of valued difference
fields with constants from $K$ which is true in $W[k]$.  As $W[k']$ is an elementary
extension of $W[k]$, this sentence must also hold in $W[k']$.  That is, the 
$R'$-special points are not Zariski dense in $Y$ contradicting our hypothesis.
\end{proof}

If we were interested only in Zariski closures computed over $K$, then we would be done now, 
though we are in a position to finish the proof of~theorem~\ref{transmain}.

\begin{proof}
Let $Y \subseteq \cX_{K'}$ be an irreducible subvariety of $\cX_{K'}$ containing 
a Zariski dense set of $R'$-special points.  
By Lemma 3.3 of~\cite{Sc au} there are a natural number $N$, 
a (possibly reducible) 
subvariety $Z \subseteq (\cX \times \cX^N)_K$, and an $R'$-special point 
$b \in \cX^N(R')$ so that $Z$ contains a Zariski dense set of $R'$-special
points and $Y = Z_b$.  By Lemma~\ref{firsttransfer}, $Z$ is a special subvariety
of $\cX^{1+N}$.  Of course, the variety $\cX \times \{ b \}$ is an $R'$-special 
subvariety of $\cX^{1+N}$ and (the components of) the intersection of two 
weakly special varieties are again weakly special.  Thus, as $Y$ contains an 
$R'$-special point, it is $R'$-special.
\end{proof}

\subsection{Bounds}

On very general grounds, Theorem~\ref{main} 
implies uniform bounds.  
At the qualitative level, we have the following uniformity.

\begin{theorem}
\label{boundqual}
Fix $g$, $n$ and $R$ as in Theorem~\ref{main} and 
a quasiprojective embedding of $\cX$.
Then there is a function $f:\NN \to \NN$ so that
for $Y \subseteq \cX_\CC$ a subvariety of degree $d$, the 
degree of the Zariski closure of the set of $R$-special points on $Y$ is 
at most $f(d)$. 
\end{theorem}
\begin{proof}
By Theorem~\ref{main} applied to the Cartesian powers of $\cX^N$, 
every irreducible subvariety of $\cX_\CC^N$ containing a Zariski
dense set of $R$-special points is $R$-special.  If $C$ is a component
of the intersection of two $R$-special varieties and $C$ contains 
an $R$-special point, then $C$ is $R$-special.  Thus, by 
Theorem 2.4 of~\cite{Sc au} the set of $R$-special points on $\cX$ satisfies
automatic uniformity and the function $f$ may be computed using Chow 
varieties.
\end{proof}

In the course of proving Theorem~\ref{main} we did not 
completely describe the sets of form $Y(\CC) \cap \Psi_m$ where
$Y \subseteq \cX_\CC$ is a closed subvariety.   We suggest the following
conjecture.

\begin{conj}
\label{weakspecialconj}
If $Y \subseteq \cX_\CC$ is an irreducible subvariety and
$Y(\CC) \cap \Psi_m$ is Zariski dense in $Y$, then $Y$ is weakly 
special.
\end{conj}

If Conjecture~\ref{weakspecialconj} holds, then 
our method of proof for Theorem~\ref{main}
provides a method for computing the function 
$f$ of Theorem~\ref{boundqual}.  

\begin{proposition}
\label{boundexplicit}
Assume Conjecture~\ref{weakspecialconj} holds.
Fix $g$, $n$, and $R$ as in Theorem~\ref{main} and 
a quasiprojective embedding of $\cX$.  Let
$\delta$ be the degree of $\widetilde{T}_p$ (with respect to this
particular embedding).  If $Y \subseteq \cX_\CC$ is a subvariety of 
degree $\deg(Y)$, then the degree of the Zariski closure of the set of
$R$-special points on $Y$ is at most 
$[(\deg(Y))^2\delta]^{2^{\min \{ \frac{g^3 + g}{2}, 2 \dim(Y) \}}}$.
\end{proposition}

\begin{proof} 
Using the notation of our proof Theorem~\ref{main}, by 
Proposition 1.2.2 of~\cite{Hr MM}, the degree of $Z$, the Zariski 
closure of $Y(\CC) \cap \Psi_1$, is at most   
$(\deg(\widetilde{T}_p) \deg(Y)^2)^{2^{\min \{ \dim (\widetilde{T}_p), 
	2 \dim(Y) \}}}$, which is the bound stated in the Proposition.
If Conjecture~\ref{weakspecialconj} holds, then every component of 
the Zariski closure, $Z'$, of the set of $R$-special points on $Y$ is 
a component of $Z$.  Hence, $\deg(Z') \leq \deg(Z)$ and the
desired inequality holds.
\end{proof}

Since we have captured the $R$-special points by difference equations
involving the relative Frobenius, there are uniform
bounds on the $p$-adic proximity from $R$-special points to subvarieties
of $X$.   We recall the definition of proximity to a variety.

\begin{definition}
Let $(K,v)$ be a valued field.
Let $X \subseteq {\mathbb A}^N_K$ be an affine variety over $K$
with a fixed choice $F_1, \ldots, F_\ell$ of generators of its 
ideal of definition.  For a point $a \in {\mathbb A}_K^N$ we define 
the proximity of $a$ to $X$ to be $\lambda_v(X,a) := 
\max\{ 0, \min \{ v(F_i(a)) \} \}$.

If $Y$ is a finitely presented algebraic variety over $K$
with a finite affine covering $Y = \bigcup U_i$  and 
$X \subseteq Y$ is a closed subvariety, then for any $a \in Y(K)$
we define $\lambda_v(X,a)$ to be the minimum value of $\lambda_v(X \cap U_i,a)$
read in the affine charts.
\end{definition}

The definition of proximity depends on the choice of generators
of $I(X)$ and the covering of $Y$, but these dependencies are irrelevant 
to the following statements.  If the variety $X$ has a good integral model, then 
it would be reasonable to fix integral generators for the ideal of definition of
$X$ in the definition of the proximity to $X$ function.  In so doing, for 
an integral point $a$ of the ambient space, $\lambda(X,a) = \gamma$ just in 
case the reduction modulo the ideal $I_\delta := \{ y \in \cO_K ~|~ v(y) > \delta \}$
lies in the reduction of $X$ if and only if $\delta < \gamma$.  

\begin{proposition}
\label{TVR}
With the notation as in Theorem~\ref{main}, if $Z \subseteq \cX_K$ is a 
closed subvariety of $\cX_K$, then there is a constant $C \in \ZZ$ 
that if $\xi \in \cX(R)$ is an $R$-special point, then either
$\xi \in Z(R)$ or $\lambda(Z,\xi) \leq C$.
\end{proposition}

\begin{proof}
This proof follows the general outline of the proof of Theorem 0.3 of~\cite{Sc tv}.
This proof is composed of two claims.

\noindent
{\bf Claim 1}:  If $Y$ is the Zariski closure of the set of $R$-special 
points on $Z$, then there is a number $N \in \NN$ so that for any 
$R$-special point $\xi \in \cX(R)$ one has $\lambda(Z,\xi) \leq N ( \lambda(Y,\xi) + 1)$.

\noindent
{\bf Proof of Claim 1:}   Suppose that the claim fails. For each 
natural number $N$ let $\xi_N$ be an $R$-special point with 
$\lambda(Z,\xi_N) > N(\lambda(Y,\xi_N) + 1)$.   Let 
$R^*$ be a nonprincipal ultrapower of $R$ and $\widetilde{\xi} \in \cX(R^*)$
the image of the infinite tuple $\langle \xi_N \rangle_{N \in \NN}$ under
the quotient map $\prod_{N \in \NN} R \to R^*$.   

We continue to denote the (Krull) valuation on the fraction field
of $R^*$ by $v$.  Let 
$I := \{y \in R^* ~|~  (\forall N \in \NN) \ v(y) > 
N (\lambda(Y,\widetilde{\xi}) + 1) \}$.  Since each $\xi_i$ is integral, 
there is an integer $M$ such that $\lambda(Y,\widetilde{\xi}) > M$.  Thus, 
$I$ is a non-unit ideal of $R^*$ and under the usual diagonal 
embedding of $R$ into the ultrapower, $I \cap R$ is the zero ideal. Set 
$\overline{R} := R^*/I$ and let $\overline{\xi}$ be the image of 
$\widetilde{\xi}$ under the quotient map.  By the above observation,
the natural map from $R$ to $\overline{R}$ is an embedding.

 The formula $\phi$
of Lemma~\ref{definable} is closed under the specialization 
$R^* \to \overline{R}$ so that $\overline{\xi}$ is $\overline{R}$-special.
By the main theorems of~\cite{Sc} or~\cite{BMS}, the extension 
$R \hookrightarrow \overline{R}$ is elementary.  As every $R$-special point
in $Z$ is contained in $Y$, the same must be true of $\overline{R}$-special
points, but by construction $\overline{\xi} \in Z(\overline{R}) \smallsetminus Y(\overline{R})$.
\hspace{\fill} $\Box$
  
\medskip
\noindent
{\bf Claim 2}: If $Y \subseteq \cX$ is an $R$-special subvariety (not necessarily irreducible) 
and $\xi \in \cX(R)$ is an $R$-special point not on $Y$, then, provided that
we compute proximity to $Y$ with respect to generators of its ideal over $R$, $\lambda(Y,\xi) \leq \frac{1}{p -1}$.

\medskip
\noindent
{\bf Proof of Claim 2:}  As the proximity to a finite union is the minimum of the proximities to each set, we may 
assume that $Y$ is irreducible.  If $B \leq \cX_Y$ is the stabilizer of $Y$, and $g:Y \to Y/B \subseteq \cX'$
is the quotient map, then $\lambda(Y,\xi) \leq \lambda(Y/B,\pi(\xi))$.  So, it suffices to consider the case that
the fibres of $Y$ over $\pi(Y)$ are finite, hence, (because $Y$ is $R$-special) torsion. 

Let $L$ be an algebraic closure of $R[\frac{1}{p}]$ and continue to denote by $v$ the extension 
of the valuation to $L$.   Write $\cO_L := \{ x \in K ~|~ v(x) \geq 0 \}$ for the 
ring of integers in $K$ and set $J := \{ x \in K ~|~ v(x) > \frac{1}{p-1} \}$.
Let $\nu:\cO_K \to \cO_K/J$ be the quotient map and continue to denote by $\nu$ the 
reduction maps $\nu:\cX(\cO_K) \to \cX(\cO_K/I)$ and $\nu:\cA(\cO_K) \to \cA(\cO_K/I)$.
Then each such $\nu$
is injective on the set of $R$-special points.   In fact, if 
$\Xi := \{ \zeta \in \cX(\cO_K) ~|~  \cX_{\pi(\zeta)} \text{ is a canonical lift and } \zeta \text{ is torsion } \}$, 
then $\nu$ is injective on $\Xi$.  Thus, if $\nu(\xi) \in Y(\cO_K/J)$, then $\xi \in Y(R)$. 
\hspace{\fill} $\Box$

\medskip 
The result follows by combining the two claims, taking $C = 2N$.

\end{proof}

In the course of proving Claim 2 we actually showed that for a special subvariety $Y \subseteq \cX_R$ 
the uniform bound on the proximity to $Y$ is valid over the set of points (over $\QQ_p^\alg$) which are 
torsion points on canonical lift fibres.  We conjecture that the restriction to $Y$ special is unnecessary.

\begin{conj}
\label{canonicalconj}
Let $\Cp$ be the completion of the algebraic closure of $\Qp$ 
with $v$ the extension of the $p$-adic valuation to $\Cp$
and $\cO = \cO_\Cp := \{ x \in \Cp ~|~ v(x) \geq 0 \}$ the ring of integers in $\Cp$. 
Let $Z \subseteq \cX_\Cp$ be a closed subvariety of $\cX_\Cp$.
Let $\lambda(Z,\cdot)$ be a $p$-adic proximity function for $Z$.  
Then there is a rational number $q \in \QQ$ so that for any
point $\xi \in \cX(\Cp)$ which is a torsion point of a canonical lift
fibre either $\xi \in Z(\Cp)$ or $\lambda(Z,\xi) \leq q$.
\end{conj}

There are many natural strengthenings of Conjecture~\ref{canonicalconj}.  
For instance, one might require of $\xi$ only that it is an integral 
special point (thereby including quasi-canonical fibres) or one might 
ask for the bound to hold for \emph{all} special points.   While this 
last formulation may seem particularly dangerous since one would have to deal
with families of abelian varieties having bad reduction, it should be noted 
that if $Z \subseteq (\cX_a)_\Cp$ where $a \in \cA(\QQ_p^\alg)$, then, in fact, such
a bound exists (see the Main Theorem of~\cite{Sc tv2}).


\begin{thebibliography}{99}

\bibitem{Andre} {\sc Y. Andr\'{e}}, {\bf $G$-Functions and geometry}, 
Aspects of Mathematics, {\bf E13}, Friedr. Vieweg \& Sohn, Braunschweig, 1989. xii+229 pp.

\bibitem{Andre curve} {\sc Y. Andr\'{e}}, 
Finitude des couples d'invariants modulaires singuliers sur une courbe alg\'{e}brique plane non modulaire,
\emph{J. Reine Angew. Math.} {\bf 505} (1998), 203--208.

\bibitem{BMS} {\sc L. B\'{e}lair}, {\sc A. Macintyre}, and {\sc T. Scanlon}, Model theory of Frobenius
on Witt vectors, preprint, 2002, \\ available at {\tt http://www.math.berkeley.edu/$\sim$scanlon/papers/papers.html}.

\bibitem{ChHr} {\sc Z. Chatzidakis} and {\sc E. Hrushovski},  Model theory of difference fields,
 \emph{Trans. Amer. Math. Soc.} {\bf 351} (1999), no. 8, 2997--3071.

\bibitem{clark} {\sc P. Clark}, Bounds for torsion on abelian varieties with integral moduli, 
preprint, 2004, {\tt arXiv:math.NT/0407264}.

\bibitem{Cohn} {\sc R. M. Cohn}, {\bf Difference algebra},  
Interscience Publishers John Wiley \& Sons, New York-London-Sydeny 1965 xiv+355 pp.

\bibitem{Edix curve} {\sc B. Edixhoven}, Special points on the product of two modular curves,
\emph{Compositio Math.} {\bf 114} (1998), no. 3, 315--328.

\bibitem{Edix hilb} {\sc B. Edixhoven}, On the Andr\'{e}-Oort conjecture for Hilbert modular surfaces,
{\bf  Moduli of abelian varieties} (Texel Island, 1999), 133--155, Progr. Math., {\bf 195}, Birkh\"{a}user, Basel, 2001.

\bibitem{EdYa} {\sc B. Edixhoven} and {\sc A. Yafaev},  Subvarieties of Shimura varieties,
\emph{Ann. of Math. (2)} {\bf 157} (2003), no. 2, 621--645.

\bibitem{Hr} {\sc E. Hrushovski}, Proof of Manin's theorem by reduction to positive characteristic,
in {\bf Model Theory and Algebraic Geometry: An introduction to E. Hrushovski's proof of the 
geometric Mordell-Lang conjecture}, (Elisabeth Bouscaren, ed.), Springer {\bf LNM 1696}, Berlin 1998, 
197 -- 205.

\bibitem{Hr MM} {\sc E. Hrushovski}, The Manin-Mumford conjecture and the model theory of difference fields,
\emph{Ann. Pure Appl. Logic} {\bf 112} (2001), no. 1, 43--115.

\bibitem{mattuck} {\sc A. Mattuck},  Abelian varieties over $p$-adic ground fields,
\emph{Ann. of Math. (2)} {\bf 62}, (1955). 92--119.

\bibitem{Messing} {\sc W. Messing}, {\bf The Crystals Associated to Barsotti-Tate Groups: with Applications
to Abelian Schemes}, Lecture Notes in Mathematics {\bf 264}, Springer-Verlag, Berlin, 1972. 

\bibitem{Milne} {\sc J. S. Milne}, Abelian varieties, in {\bf Arithmetic Geometry}, (Gary Cornell and 
Joseph H. Silverman, eds.), Springer-Verlag, New York, 1986, 103 -- 150.

\bibitem{Moonen} {\sc B. Moonen}, Linearity properties of Shimura varieties, part II, 
\emph{Compositio Math.} {\bf 114} (1998), no. 1, 3--35.


\bibitem{GIT} {\sc D. Mumford}, {\sc J. Fogarty}, and {\sc F. Kirwan},
 {\bf Geometric Invariant Theory}, third edition, Ergebnisse der Mathematik
und ihrer Grenzgebiete {\bf 34}, Spinger-Verlag (New York), 1994.

\bibitem{Oort} {\sc F. Oort}, Canonical liftings and dense sets of CM-points, in  
{\bf Arithmetic geometry} (Cortona, 1994), 228--234, 
Sympos. Math., XXXVII, Cambridge Univ. Press, Cambridge, 1997.

\bibitem{Pink} {\sc R. Pink}, {\bf Arithmetical Compactifications of Mixed Shimura Varieties}, 
Ph.D. thesis, Friedrich-Wilhelms-Universit\"{a}t, Bonner Mathematische Schriften {\bf 209}, 1989.

\bibitem{Pink aoml} {\sc R. Pink}, A Combination of the Conjectures of Mordell-Lang and Andr\'{e}-Oort,
preprint, August 2004, available at \\ {\tt http://www.math.ethz.ch/$\sim$pink/preprints.html}

\bibitem{PiRo} {\sc R. Pink} and {\sc D. Roessler}, On Hrushovski's proof of the Manin-Mumford conjecture,
Proceedings of the International Congress of Mathematicians, Vol. I (Beijing, 2002), 539--546, 
Higher Ed. Press, Beijing, 2002.

\bibitem{Sc tv} {\sc T. Scanlon}, $p$-adic distance from torsion points of semi-abelian varieties, 
\emph{J. Reine Angew. Math.} {\bf 499} (1998), 225 -- 236.

\bibitem{Sc tv2} {\sc T. Scanlon}, The conjecture of Tate and Voloch on $p$-adic proximity to torsion,
 \emph{International Mathematics Research Notices} {\bf 1999}, no. 17, 909 -- 914. 

\bibitem{Sc} {\sc T. Scanlon}, Quantifier elimination for the relative Frobenius, in 
{\bf Valuation Theory and Its Applications Volume II},
conference proceedings of the International Conference on Valuation Theory (Saskatoon, 1999),
 Franz-Viktor Kuhlmann, Salma Kuhlmann, and Murray Marshall, eds., Fields Institute Communications Series, 
(AMS, Providence), 2003, 323 - 352.

\bibitem{Sc au} {\sc T. Scanlon}, Automatic uniformity, \emph{International Mathematics Research Notices},
(to appear).

\end{thebibliography}
\end{document}